\documentclass[12pt,leqno]{article}

\usepackage{amssymb}
\usepackage{latexsym}
\usepackage{exscale}

\title{Quantized orthonormal systems: \\ A
non-commutative Kwapie\'n theorem}

\author{{\sc J. Garc\'{\i}a-Cuerva} and {\sc J.Parcet}}

\date{}

\newtheorem{theorem}{Theorem}[section]
\newtheorem{lemma}[theorem]{Lemma}
\newtheorem{definition}[theorem]{Definition}
\newtheorem{corollary}[theorem]{Corollary}
\newtheorem{remark}[theorem]{Remark}
\newtheorem{proposition}[theorem]{Proposition}

\newcommand{\N}{\mathbb{N}}
\newcommand{\Z}{\mathbb{Z}}
\newcommand{\C}{\mathbb{C}}

\newcommand{\sig}{\sigma}
\newcommand{\Sig}{\Sigma}
\newcommand{\Degree}{d_{\sigma}}

\newcommand{\Fourier}{\mathcal{F}}

\newcommand{\Lebesgue}{\mathcal{L}}

\def\essup{\mathop{\rm ess\, sup}}

\begin{document}

\textwidth=13.5cm

\maketitle

 \footnote{Research supported in part by the
European Commission via the TMR Network `Harmonic Analysis' and by
Project BFM 2001/0189, Spain}

\begin{abstract}
The concepts of Riesz type and cotype of a given Banach space are
extended to a non-commutative setting. First, the Banach space is
replaced by an ope\-rator space. The notion of quantized
orthonormal system, which plays the role of the orthonormal system
in the classical setting, is then defined. The Fourier type and
cotype of an ope\-rator space with respect to a non-commutative
compact group fit in this context. Also, the quantized analogs of
Rademacher and Gaussian systems are treated. All this is used to
obtain an operator space version of the classical theorem of
Kwapie\'n characterizing Hilbert spaces by means of vector-valued
orthogonal series. Several approaches to this result with
different consequences are given.
\end{abstract}

\section{Introduction}
\label{Section-Introduction}

\large\normalsize

The notion of type or cotype of a Banach space $B$ with respect to
some classical system, such as the Rademacher or the trigonometric
system, is a common way to express the validity of certain
inequalities for $B$-valued functions. The systematic research of
these topics has given rise to a very well developed theory of the
interaction between orthonormal systems and the geometry of Banach
spaces. In this paper we look at this interaction from a
non-commutative point of view. By that we mean to investigate what
happens when we replace Banach spaces by operator spaces.

The first example in this setting was given in \cite{GP}, where we
define and study the Fourier type and cotype of an operator space
with respect to a non-commutative compact group. Let $1 \le p \le
2$ and let $p'$ denote its conjugate exponent. Let $G$ be a
compact group with dual object $\Gamma$. An operator space $E$ has
Fourier type $p$ with respect to $G$ if the $E$-valued Fourier
transform on $G$ extends to a completely bounded operator from
$L_E^p(G)$ into $\Lebesgue_E^{p'}(\Gamma)$. Similarly, by
considering the inverse of the Fourier transform, the notion of
Fourier cotype comes out in this context. A relevant difference,
between this notion of Fourier type and its classical counterpart
for compact abelian groups, lies in the fact that the system of
characters has to be replaced by the set of equivalence classes of
unitary irreducible representations of $G$. That is, $\Gamma$ is
now the system and $E$ is the space.

Going back to the general case, the question is to find out which
properties should we require from the system to get appropriate
information about the operator space. As we shall see below, these
systems will be collections of matrix-valued functions satisfying
some extra conditions, what is perfectly natural in view of the
basic example mentioned above. This is why we have called them \lq
quantized systems\rq, completing in such a way the scheme where
Banach spaces become operator spaces and boundedness of operators
is replaced by complete boundedness. Finally we point out that, by
the necessity of working with vector-valued Schatten classes and
as it was recalled in \cite{P2}, the management of vector-valued
orthogonal series with respect to a quantized system does not make
sense in Banach space theory.

The definition of quantized orthonormal system was motivated by
the theory initiated in \cite{GP} and \cite{GMP}. If fact, we find
basic for its development to obtain an operator space version of
the isomorphic characterization of Hilbert spaces given by
Kwapie\'n in \cite{K1}. We provide three different approaches to
this result. The first one is valid for any uniformly bounded
quantized orthonormal system. The second one extends this result
to non-uniformly bounded but complete quantized orthonormal
systems. The third approach involves the quantization of the
classical Gauss system, which fails to be complete or even
uniformly bounded. This system also characterizes Pisier's OH
Hilbertian operator spaces up to complete isomorphism and the
proof of this fact follows the arguments given in the first
approach. However, we also show that Kwapie\'n's original
arguments in \cite{K1} to link Rademacher and Gauss systems via
the central limit theorem work in this context. Moreover, as we
shall see, the use of this probabilistic approach has a remarkable
advantage. Namely, it provides corollary \ref{Kwapien-Banach}.
Roughly speaking this result can be stated by saying that, when
the quantized system we deal with takes values in arbitrary large
matrices, then the operator space version of Kwapie\'n's theorem
for such a system also holds requiring only the boundedness of the
involved operators, not the complete boundedness. Finally, an
example is included and some open questions are posed.

All throughout this paper some basic notions of operator space
theory and vector-valued Schatten classes will be assumed. The
definitions and results about operator spaces that we will be
using can be found in the book of Effros and Ruan \cite{ER}, while
for the study of those Schatten classes the reader is referred to
\cite{P2}, where Pisier analyzes them in detail.

\

\textbf{Acknowledgment.} We wish to thank Gilles Pisier for some
valuable comments. Specially for the proof of theorem
\ref{Complete}, which he communicated to us during a visit of the
second author to Universit\'{e} Paris VI.

\section{Uniformly bounded quantized orthonormal systems}
\label{Section-Definitions}

The classical Hausdorff-Young inequality on the torus was
generalized by F. Riesz in $1923$ to any uniformly bounded
orthonormal system. If one looks for extensions of this result to
vector-valued functions, the notions of Riesz type, cotype and
strong cotype of a Banach space come out naturally. These were
defined in \cite{GKKT} with the aim to provide a general notion of
type which included the classical (uniformly bounded) systems:
Rademacher, Fourier, Walsh, etc... Here we introduce a \lq
quantized version\rq${}$ of these notions. From now on, $M_n$ will
stand for the vector space of $n \times n$ complex matrices and
$S_n^p$ will denote the Schatten $p$-class over the space $M_n$.

\begin{definition} \label{System}
\emph{Let $(\Omega,\mathcal{M},\mu)$ be a probability measure
space with no atoms and let $\mathbf{d}_{\Sig} = \{\Degree: \sig
\in \Sig\}$ be a family of positive integers, $\Sig$ an index set.
A collection of matrix-valued functions $\Phi = \{\varphi^{\sig}:
\Omega \rightarrow M_{\Degree}\}_{\sig \in \Sig}$ with
measura\-ble entries is said to be a \emph{uniformly bounded
quantized orthonormal system} (u.b.q.o.s. for short) if the
following conditions hold:
\begin{itemize}
\item[{\bf (a)}] $\ \ \displaystyle \int_{\Omega}
\varphi_{ij}^{\sig} (\omega) \overline{\varphi_{i'j'}^{\sig'}
(\omega)} d \mu(\omega) = \,\ \frac{1}{\Degree}  \,\
\delta_{\sig\sig'}\delta_{ii'}\delta_{jj'}.$ \\ ${}$
\item[{\bf (b)}] $\ \ \displaystyle \sup_{\sig \in \Sig} \,\
\essup_{\omega \in \Omega} \,\
\|\varphi^{\sig}(\omega)\|_{S_{\Degree}^{\infty}} = \,\
\mathrm{M}_{\Phi} < \infty$.
\end{itemize}
The pair $(\Sig, \mathbf{d}_{\Sig})$ will be called the \emph{set
of parameters} of $\Phi$. We say that $\Phi$ is \emph{complete}
when any function $f \in L^2(\Omega)$ can be written as $$f =
\sum_{\sig \in \Sig} \Degree \textnormal{tr}(A^{\sig}
\varphi^{\sig}) \ \ \mbox{for some} \ \ A \in \prod_{\sig \in
\Sig} M_{\Degree}.$$}
\end{definition}

\begin{remark}
\emph{Let us recall that, if we take $\Sig = \N$ and $\Degree = 1$
for all $\sig \in \Sig$, we recover the classical definition of
uniformly bounded orthonormal systems or complete orthonormal
systems on $\Omega$. Also, if $\Omega$ is a compact topological
group $G$ with normalized Haar measure $\mu$, then the dual object
$\Gamma$ of $G$ is a u.b.q.o.s. The functions $\varphi^{\sig}$ are
irreducible unitary representations of $G$, $\Degree$ is the
degree of $\varphi^{\sig}$ and $\mathrm{M}_{\Gamma} = 1$.}
\end{remark}

Let $1 \le p < \infty$, let $E$ be an operator space and let
$\Sig$ be an index set as in definition \ref{System}. Following
the notation of \cite{GP} we define the spaces $$\begin{array}{l}
\displaystyle \Lebesgue_E^p(\Sig) = \Big\{ A \in \prod_{\sig \in
\Sig} M_{\Degree} \otimes E: \,\ \|A\|_{\Lebesgue_E^p(\Sig)} =
\Big( \sum_{\sig \in \Sig} \Degree
\|A^{\sig}\|_{S_{\Degree}^p(E)}^p \Big)^{1/p} < \infty \Big\} \\
\displaystyle \Lebesgue_E^{\infty}(\Sig) = \Big\{ A \in
\prod_{\sig \in \Sig} M_{\Degree} \otimes E: \,\
\|A\|_{\Lebesgue_E^{\infty}(\Sig)} = \,\ \sup_{\sig \in \Sig} \,\
\|A^{\sig}\|_{S_{\Degree}^{\infty}(E)} < \infty \Big\}
\end{array}$$ where we write $S_n^p(E)$ for the
$E$-valued Schatten $p$-class over $M_n$. $\Lebesgue^p(\Sig)$ will
denote the scalar-valued case. $\Lebesgue_E^p(\Sig)$ is endowed
with its natural operator space structure, see \cite{GP} and
Chapter $2$ of \cite{P2} for the details. Now, if $\Phi$ is a
u.b.q.o.s. and $\star$ stands for the adjoint operator, then the
\emph{$\Phi$-transform} and its inverse can be defined naturally
as follows: $$\Fourier_{\Phi}(f)^{\sig} = \int_{\Omega} f(\omega)
\varphi^{\sig}(\omega)^{\star} d \mu(\omega) \quad \mbox{and}
\quad \Fourier_{\Phi}^{-1}(A) (\omega) = \sum_{\sig \in \Sig}
\Degree \mbox{tr}(A^{\sig} \varphi^{\sig}(\omega))$$ for $f:
\Omega \rightarrow E$ and $\displaystyle A \in \prod_{\sig \in
\Sig} M_{\Degree} \otimes E$.

We start with a version for uniformly bounded quantized
orthonormal systems of the classical Riesz theorem.

\begin{lemma} \label{Scalar}
Let $1 \le p \le 2$ and let $p'$ denote its conjugate exponent.
Let $\Phi$ be a u.b.q.o.s. Then we have
$$\|\Fourier_{\Phi}\|_{cb(L^p(\Omega), \Lebesgue^{p'}(\Sig))}, \,\
\|\Fourier_{\Phi}^{-1}\|_{cb(\Lebesgue^p(\Sig), L^{p'}(\Omega))}
\le \,\ \mathrm{M}_{\Phi}^{2/p - 1}.$$
\end{lemma}

\emph{Proof}. By the complex interpolation method for operator
spaces, we just need to check the cases $p = 1,2$. It follows from
Lemma $1.7$ of \cite{P2} that
$$\|\Fourier_{\Phi}\|_{cb(L^2(\Omega), \Lebesgue^2(\Sig))} =
\sup_{n \ge 1} \|\Fourier_{\Phi} \otimes
I_{M_n}\|_{\mathcal{B}(L_{S_n^2}^2(\Omega),
\Lebesgue_{S_n^2}^2(\Sig))}$$ with the obvious modifications for
the inverse operator. Then the case $p = 2$ is a consequence of
the orthonormality of $\Phi$. That is, it follows from condition
$(\mathbf{a})$ in definition \ref{System}. If $p = 1$,
$\Fourier_{\Phi}$ is defined on $L^1(\Omega)$ which is equipped
with the max quantization. Moreover, $\Fourier_{\Phi}^{-1}$ takes
values in $L^{\infty}(\Omega)$, which is equipped with the min
quantization. Therefore, by the quantizations we are working with,
boundedness is equivalent to complete boundedness (see Section
$3.3$ of \cite{ER} for the details). But it is obvious that the
stated inequalities hold for $p = 1$ when the $cb$ norm is
replaced by the operator norm. $\blacksquare$

\

If $\Sig_0$ is a finite subset of $\Sig$, let $\Phi_E^p(\Sig_0) =
\mbox{span} \{\varphi_{ij}^{\sig}: \sig \in \Sig_0\} \otimes E$
regarded as a subspace of $L_E^p(\Omega)$ with its natural
operator space structure. Also, let $\Phi_0 = \{\varphi^{\sig}:
\Omega \rightarrow M_{\Degree}\}_{\sig \in \Sig_0}$ be the
restriction of $\Phi$ to $\Sig_0$. Then $\Phi_0$ is also a
u.b.q.o.s. and lemma \ref{Scalar} holds for $\Phi_0$.

\begin{definition} \label{Riesz}
\emph{Let $1 \le p \le 2$ and let $p'$ denote its conjugate
exponent:
\begin{itemize}
\item The operator space $E$ is said to have \emph{Riesz type
$p$} with respect to $\Phi$, or simply $\Phi$-type $p$, if
$$\mathcal{K}_{1p}(E,\Phi) = \sup \|\Fourier_{\Phi_0}^{-1} \otimes
I_E \|_{cb(\Lebesgue_E^p(\Sig_0),\Phi_E^{p'}(\Sig_0))} < \infty$$
where the supremum is taken over the family of finite subsets
$\Sig_0$ of $\Sig$.
\item The operator space $E$ is said to have \emph{Riesz cotype $p'$}
with respect to $\Phi$, or simply $\Phi$-cotype $p'$, if
$$\mathcal{K}_{2p'}(E,\Phi) = \sup \|\Fourier_{\Phi_0} \otimes I_E
\|_{cb(\Phi_E^p(\Sig_0),\Lebesgue_E^{p'}(\Sig_0))} < \infty.$$ The
supremum is taken again over the family of finite subsets $\Sig_0$
of $\Sig$.
\item The operator space $E$ is said to have \emph{strong Riesz
cotype $p'$} with respect to $\Phi$, or simply strong
$\Phi$-cotype $p'$, if $$\mathcal{K}_{3p'}(E,\Phi) =
\|\Fourier_{\Phi,E}\|_{cb(L_E^p(\Omega),\Lebesgue_E^{p'}(\Sig))} <
\infty$$ where $\Fourier_{\Phi,E}$ denotes the extension of
$\Fourier_{\Phi} \otimes I_E$ to $L_E^p(\Omega)$.
\end{itemize}}
\end{definition}

\begin{remark}
\emph{Let us note that, if $E$ has $\Phi$-type $p$, then in
particular there exists a positive constant $c$ such that $$\Big\|
\sum_{\sig \in \Sig_0} \Degree \mbox{tr}(A^{\sig} \varphi^{\sig})
\Big\|_{L_E^{p'}(\Omega)} \le c \,\ \Big( \sum_{\sig \in \Sig_0}
\Degree \|A^{\sig}\|_{S_{\Degree}^p(E)}^p \Big)^{1/p}$$ for any
finite subset $\Sig_0$ of $\Sig$ and any $A \in
\Lebesgue_E^p(\Sig_0)$. This expression is now much closer to the
classical notion of Riesz type. In fact, for $\Degree = 1$ and
$\Sig_0 = \{1,2, \ldots n\}$, we recover the classical definition.
Analogous remarks hold for the Riesz cotype and the strong Riesz
cotype.}
\end{remark}

\begin{remark}
\emph{We point out here that, as in the classical theory, a notion
of strong Riesz type would be superfluous since it would coincide
with that of Riesz type. The proof of this fact is an easy
consequence of the density of the subspace of
$\Lebesgue_E^p(\Sig)$ formed by the elements $A$ with finite
support, that is with $A^{\sigma} \neq 0$ only for finitely many
$\sig \in \Sig$. In fact,
$$\mathcal{K}_{1p}(E, \Phi) =
\|\Fourier_{\Phi,E}^{-1}\|_{cb(\Lebesgue_E^p(\Sig),
L_E^{p'}(\Omega))}.$$ Again as in the classical case, this
equivalence is not necessarily valid for the cotype. Moreover, we
have the obvious estimate $\mathcal{K}_{2p'}(E,\Phi) \le
\mathcal{K}_{3p'}(E,\Phi)$ for any u.b.q.o.s. $\Phi$, any operator
space $E$ and any $1 \le p \le 2$. However, the $\Phi$-cotype is
equivalent to the strong $\Phi$-cotype when $\Phi$ is complete. In
this paper we shall mainly be concerned with the Riesz type and
cotype. We have defined the strong Riesz cotype because, as we
shall see below, it is the right notion for duality.}
\end{remark}

\begin{remark}
\emph{Sometimes in the sequel we shall also use the notion of
$\Psi$-type $2$ and $\Psi$-cotype $2$ for some quantized
orthonormal systems $\Psi$ which fail to be uniformly bounded.}
\end{remark}

These definitions are illustrated in \cite{GP} and \cite{GMP}
where the Fourier type and cotype of an operator space with
respect to a compact group are investigated. Namely, if $G$ is a
compact group with dual object $\Gamma$, then Fourier type $p$
with respect to $G$ is nothing but the $\Gamma$-cotype $p'$ (or
strong $\Gamma$-cotype $p'$ since $\Gamma$ is complete in $L^2(G)$
by the Peter-Weyl theorem). Moreover, Fourier cotype $p'$ with
respect to $G$ coincides with $\Gamma$-type $p$. This conflict in
our terminology goes back to the commutative theory, where Fourier
type $p$ with respect to the torus $\mathbb{T}$ means $\Z$-cotype
$p'$ (or equivalently strong $\Z$-cotype $p'$), see \cite{GKKT}.

In what follows we assume the reader is familiar with the
properties of Fourier type and cotype stated in \cite{GP} and
\cite{GMP}. In fact, we omit the proof of the following results,
since the arguments to be used can be found there.

\begin{itemize}
\item[$(a)$] \textbf{Trivial exponents}. Every operator space has
Riesz type $1$ and strong Riesz cotype $\infty$ with respect to
any u.b.q.o.s. $\Phi$. Moreover, we have the estimates
$\mathcal{K}_{11}(E,\Phi), \,\ \mathcal{K}_{2\infty}(E,\Phi), \,\
\mathcal{K}_{3\infty}(E,\Phi) \le \mathrm{M}_{\Phi}$.
\item[$(b)$] \textbf{Subspaces}. The Riesz type is preserved when
passing to subspaces. Moreover, $\mathcal{K}_{1p}(F,\Phi) \le
\mathcal{K}_{1p}(E,\Phi)$ for any closed subspace $F$ of $E$. The
same holds for the Riesz cotype and the strong Riesz cotype.
\item[$(c)$] \textbf{Complex interpolation}. Let $0 < \theta < 1$
and let $E_0$ and $E_1$ be opera\-tor spaces having $\Phi$-type
$p_0$ and $p_1$ respectively. Then, if $(E_0, E_1)$ is compatible
for complex interpolation, the interpolated operator space
$(E_0,E_1)_{\theta}$ has $\Phi$-type $p_{\theta} = p_0p_1 ((1 -
\theta) p_1 + \theta p_0)^{-1}$. In particu\-lar, the Riesz type
$p$ becomes a stronger condition on a given operator space as $p$
approaches $2$. Similar assertions hold for the Riesz cotype and
the strong Riesz cotype.
\item[$(d)$] \textbf{Duality}. $\mathcal{K}_{1p}(E,\Phi) =
\mathcal{K}_{3p'}(E^{\star},\Phi)$ and
$\mathcal{K}_{1p}(E^{\star},\Phi) = \mathcal{K}_{3p'}(E,\Phi)$.
That is, Riesz type and strong Riesz cotype are dual notions.
\item[$(e)$] \textbf{Local theory}.
If $d_{cb}$ stands for the $cb$-distance between two operator
spaces, we have $\mathcal{K}_{1p}(E_2,\Phi) \le d_{cb}(E_1,E_2)
\,\ \mathcal{K}_{1p}(E_1,\Phi)$. The same holds for the Riesz
cotype and the strong Riesz cotype.
\item[$(f)$] \textbf{Degenerate case}. Let us assume that the
index set $\Sig$ associated to the u.b.q.o.s. $\Phi$ is finite,
then we have $$\mathcal{K}_{1p}(E,\Phi), \,\
\mathcal{K}_{2p'}(E,\Phi), \,\ \mathcal{K}_{3p'}(E,\Phi) \le
\mathrm{M}_{\Phi} \,\ \Big( \sum_{\sig \in \Sig} \Degree^2
\Big)^{1/{p'}}.$$
\item[$(g)$] \textbf{Lebesgue spaces}. Let $1 \le q \le \infty$.
Let $(X,\mathcal{N},\nu)$ be a $\sigma$-finite measure space, then
$L^q(X)$ has $\Phi$-type $\min(q,q')$ and strong $\Phi$-cotype
$\max(q,q')$. Similar results hold for Schatten classes. Moreover,
$L_E^q(X)$ and $S^q(E)$ have $\Phi$-type $\min(q,q')$ and strong
$\Phi$-cotype $\max(q,q')$ whenever $E$ does.
\end{itemize}

\begin{remark}
\emph{In what follows we shall assume that $\Sig$ is not finite.}
\end{remark}

\section{The Kwapie\'n theorem for operator spaces}
\label{Section-Kwapien1}

We begin by defining the quantized version of the classical
Rademacher system. This notion is extracted from \cite{MP}, where
the authors use it to study random Fourier series on
non-commutative compact groups. From now on we fix a probability
measure space $(\Omega, \mathcal{M}, \mu)$ with no atoms, an index
set $\Sig$ and a family of positive integers $\mathbf{d}_{\Sig}$.

\begin{definition}
\emph{The quantized \emph{Rademacher} system associated to $(\Sig,
\mathbf{d}_{\Sig})$ is defined by a collection $\mathcal{R} =
\{\rho^{\sig}: \Omega \rightarrow O(\Degree)\}_{\sig \in \Sig}$ of
independent random orthogonal matrices, uniformly distributed on
the orthogonal group $O(\Degree)$ equipped with its normalized
Haar measure $\nu_{\sig}$.}
\end{definition}

\begin{remark}
\emph{Similarly, the quantized \emph{Steinhaus} system associated
to $(\Sig, \mathbf{d}_{\Sig})$ is a collection $\mathcal{S} =
\{\xi^{\sig}: \Omega \rightarrow U(\Degree)\}_{\sig \in \Sig}$ of
independent random unitary matrices, uniformly distributed on the
unitary group $U(\Degree)$ equipped with its normalized Haar
measure $\lambda_{\sig}$. It is easy to check that both Rademacher
and Steinhaus systems are u.b.q.o.s.'s with uniform bound
$\mathrm{M}_{\mathcal{R}} = \mathrm{M}_{\mathcal{S}} = 1$.
Moreover, the notions of $\mathcal{R}$-type $p$ and
$\mathcal{S}$-type $p$ are equivalent for $1 \le p \le 2$. Namely,
the inequalities $$\frac{1}{2} \,\
\|\Fourier_{\mathcal{R}}^{-1}(A)\|_{L_B^q(\Omega)} \le \,\
\|\Fourier_{\mathcal{S}}^{-1}(A)\|_{L_B^q(\Omega)} \le \,\ 2 \,\
\|\Fourier_{\mathcal{R}}^{-1}(A)\|_{L_B^q(\Omega)}$$ were proved
in \cite{MP} for any Banach space $B$, any $A$ supported in any
finite subset $\Sig_0$ of $\Sig$ and any $1 \le q < \infty$.
Hence, given an operator space $E$, we just need to take $B =
S_n^{p'}(E)$ for any $n \ge 1$ and $q =p'$ to see this
equiva\-lence. Similar arguments are valid to show that the same
equivalence holds between $\mathcal{R}$-cotype and
$\mathcal{S}$-cotype. Moreover, the equivalence between both
systems with respect to the strong Riesz cotype follows by
duality. Therefore, although the results obtained will be valid
for both systems, we shall work only with the quantized Rademacher
system.}
\end{remark}

\begin{remark} \label{Khintchine-Kahane}
\emph{Let $\mathrm{R}_p(E)$ be the closure in $L_E^p[0,1]$ of the
subspace of linear combinations of the classical Rademacher
functions $r_1, r_2, \ldots$ with $E$-valued coefficients. In
particular, we shall write $\mathrm{R}_p$ for the closure in
$L^p[0,1]$ of the subspace spanned by $r_1, r_2 \ldots$ The
classical Khintchine-Kahane inequalities can be rephrased by
saying that the norm of $\mathrm{R}_p(E)$, regarded as a Banach
space, is equivalent to that of $\mathrm{R}_q(E)$ whenever $1 \le
p \neq q < \infty$. In particular we can put any exponent $1 \le q
< \infty$ in the defining inequality of Rademacher type $p$ (resp.
cotype $p'$) for (the underlying Banach space of) $E$
\begin{equation} \label{Independence}
c_1 \,\ \Big( \sum_{k=1}^n \|e_k\|_E^{p'} \Big)^{1/p'} \le \Big\|
\sum_{k=1}^n e_k r_k \Big\|_{L_E^q[0,1]} \le c_2 \,\ \Big(
\sum_{k=1}^n \|e_k\|_E^p \Big)^{1/p}.
\end{equation}
On the other hand $\mathrm{R}_p(E)$ has a natural operator space
structure inherited from $L_E^p[0,1]$. It is a remarkable fact
that the norm of $\mathrm{R}_p(E)$ is not completely equivalent to
that of $\mathrm{R}_q(E)$. That is, the operator spaces
$\mathrm{R}_p(E)$ and $\mathrm{R}_q(E)$ are isomorphic but not
completely isomorphic. The proof of this fact is due to Pisier and
it can be found in Chapter $8$ of \cite{P2}. If we replace $r_1,
r_2, \ldots$ by the entries of a quantized Rademacher system
$\mathcal{R}$, then we obtain an operator space $\mathcal{R}_p(E)$
which is Banach isomorphic but not completely isomorphic to
$\mathcal{R}_q(E)$ whenever $1 \le p \neq q < \infty$. This
equivalence of the norms, which fails to be complete, follows from
a version of the Khintchine-Kahane inequalities for $\mathcal{R}$
stated in \cite{MP}. Therefore, in contrast with
(\ref{Independence}), each election of the exponent $1 \le q <
\infty$ in definition \ref{Riesz} gives different notions of
Rademacher type and cotype. For instance, one could be tempted to
take $q$ to be $2$ no matter which would be the value of $p$. In
fact, this alternative definition becomes very useful in some
other contexts which do not appear in this paper, such as the
study of the notion of \emph{non-trivial Rademacher type}. In any
case, we have no risk in this paper to choose the wrong definition
since we shall mainly be concerned with the quadratic case $p=2$.}
\end{remark}

Now we prove the extremality of the quantized Rade\-macher system
with respect to Riesz type and cotype among the family of
uniformly bounded quantized orthonormal systems. We shall need the
following version, given in \cite{MP}, of the classical
contraction principle.

\begin{lemma} \label{Contraction}
Let $B$ be a Banach space, $\Sig_0 \subset \Sig$ finite, $A^{\sig}
\in M_{\Degree} \otimes B$ and $D^{\sig} \in M_{\Degree}$ for
$\sig \in \Sig_0$. Then, for any $1 \le q < \infty$ we have
$$\Big\| \sum_{\sig \in \Sig_0} \Degree \textnormal{tr} (A^{\sig}
\rho^{\sig} D^{\sig}) \Big\|_{L_B^q(\Omega)} \le \sup_{\sig \in
\Sig_0} \|D^{\sig}\|_{S_{\Degree}^{\infty}} \Big\|\sum_{\sig \in
\Sig_0} \Degree \textnormal{tr}(A^{\sig} \rho^{\sig})
\Big\|_{L_B^q(\Omega)}.$$
\end{lemma}

\begin{proposition} \label{Extremal}
Let $1 \le p \le 2$ and let $p'$ denote its conjugate exponent.
Then, the following holds for any operator space $E$ and any
$u.b.q.o.s.$ $\Phi\!:$
\begin{itemize}
\item[$1.$] If $E$ has $\Phi$-type $p$, then $E$ has $\mathcal{R}$-type
$p$.
\item[$2.$] If $E$ has $\Phi$-cotype $p'$, then $E$ has $\mathcal{R}$-cotype
$p'$.
\end{itemize}
\end{proposition}

\emph{Proof}. The case $p = 1$ is trivial, hence we assume that
$E$ has $\Phi$-type $p$ for some $1 < p \le 2$. First we recall
the completely isometric isomorphism
\begin{equation} \label{Fubini}
S_n^{p'}(L_E^{p'}(\Omega)) = L_{S_n^{p'}(E)}^{p'}(\Omega).
\end{equation}
On the other hand, by the orthonormality of $\Phi$ we have
\begin{equation} \label{Orthonormal}
\int_{\Omega} |\varphi^{\sig}|^2 d\mu = I_{M_{\Degree}}
\end{equation}
for all $\sig \in \Sig$. Hence, given $n \ge 1$ and $A_{ij} \in
\Lebesgue_E^p(\Sig_0)$ for $1 \le i,j \le n$, we apply
(\ref{Fubini}), (\ref{Orthonormal}), Jensen's inequality and the
contraction principle stated in lemma \ref{Contraction} to get
\begin{eqnarray*} \lefteqn{\Big\| \Big( \sum_{\sig \in \Sig_0}
\Degree \mbox{tr} (A_{ij}^{\sig} \rho^{\sig}) \,\ \Big)
\Big\|_{S_n^{p'}(L_E^{p'}(\Omega))}} \\ & = & \Big[ \int_{\Omega}
\Big\| \int_{\Omega} \Big( \sum_{\sig \in \Sig_0} \Degree
\mbox{tr} \big[ (\rho^{\sig}(\omega_1)
|\varphi^{\sig}(\omega_2)|^2 A_{ij}^{\sig} \big] \,\ \Big)
d\mu(\omega_2) \Big\|_{S_n^{p'}(E)}^{p'} d\mu(\omega_1) \,\
\Big]^{1/p'} \\ & \le & \Big[ \int_{\Omega} \int_{\Omega} \Big\|
\Big( \sum_{\sig \in \Sig_0} \Degree \mbox{tr} \big[
\rho^{\sig}(\omega_1) \varphi^{\sig}(\omega_2)^{\star}
\varphi^{\sig}(\omega_2) A_{ij}^{\sig} \big] \Big)
\Big\|_{S_n^{p'}(E)}^{p'} d\mu(\omega_1) d\mu(\omega_2)
\Big]^{1/p'} \\ & \le & \mathrm{M}_{\Phi} \Big[ \int_{\Omega}
\int_{\Omega} \Big\| \Big( \sum_{\sig \in \Sig_0} \Degree
\mbox{tr} \big[ (\rho^{\sig}(\omega_1) \varphi^{\sig}(\omega_2)
A_{ij}^{\sig} \big] \,\ \Big) \Big\|_{S_n^{p'}(E)}^{p'}
d\mu(\omega_2) d\mu(\omega_1) \Big]^{1/p'} \\ & \le &
\mathrm{M}_{\Phi} \,\ \mathcal{K}_{1p}(E,\Phi) \,\ \Big[
\int_{\Omega} \Big\| \Big( \,\ A_{ij}^{\sig} \rho^{\sig}(\omega_1)
\,\ \Big)_{\sig \in \Sig_0}
\Big\|_{S_n^{p'}(\Lebesgue_E^p(\Sig_0))}^{p'} d\mu(\omega_1)
\Big]^{1/p'}
\end{eqnarray*}
Finally, by virtue of Lemma $1.7$ of \cite{P2}, it remains to see
that $$\Big[ \int_{\Omega} \Big\| \Big( \,\ A_{ij}^{\sig}
\rho^{\sig}(\omega_1) \,\ \Big)_{\sig \in \Sig_0}
\Big\|_{S_n^{p'}(\Lebesgue_E^p(\Sig_0))}^{p'} d\mu(\omega_1)
\Big]^{1/p'} = \Big\| \Big( \,\ A_{ij} \,\ \Big)
\Big\|_{S_n^{p'}(\Lebesgue_E^p(\Sig_0))}.$$ To that aim it
suffices to check that the mapping $A \mapsto (A^{\sig}
\rho^{\sig} (\omega))_{\sig \in \Sig_0}$ is a complete isometry
from $\Lebesgue_E^p(\Sig_0)$ into itself. But this follows from
the fact that $\rho^{\sig}(\omega) \in O(\Degree)$ for all $\sig
\in \Sig$ and all $\omega \in \Omega$, see Lemma $1.6$ of
\cite{P2}. This gives the estimate
$\mathcal{K}_{1p}(E,\mathcal{R}) \le \mathrm{M}_{\Phi} \,\
\mathcal{K}_{1p}(E,\Phi)$. Similar arguments give the relation
$\mathcal{K}_{2p'}(E,\mathcal{R}) \le \mathrm{M}_{\Phi} \,\
\mathcal{K}_{2p'}(E,\Phi)$. This completes the proof.
$\blacksquare$

\begin{remark}
\emph{By duality, a similar result holds for the strong cotype.}
\end{remark}

The following is a classical result which characterizes, in terms
of the convergence of some series of vector-valued random
variables, Rademacher type (resp. cotype) $2$ Banach spaces. The
proof can be found in \cite{AG}, see Theorem $7.2$ of Chapter $3$.

\begin{lemma} \label{Araujo-Gine}
The following assertions hold:
\begin{itemize}
\item[$1.$] The Banach space $B$ has Rademacher type $2$ if and
only if there exists a sequence $\zeta_1, \zeta_2, \ldots$ of mean
zero independent random variables in $L^2(\Omega)$ with $0 < c_1
\le \|\zeta_n\|_{L^2(\Omega)} \le c_2 < \infty$ and such that, if
$x_1, x_2, \ldots$ is any sequence in $B$, then we have
$$\sum_{k=1}^{\infty} \|x_k\|_B^2 < \infty \Longrightarrow
\sum_{k=1}^n x_k \zeta_k \ \ \mbox{converges} \ \
\textnormal{a.s.}$$
\item[$2.$] The Banach space $B$ has Rademacher cotype $2$ if
and only if there exists a sequence $\zeta_1, \zeta_2, \ldots$ of
mean zero independent random variables in $L^2(\Omega)$ with $0 <
c_1 \le \|\zeta_n\|_{L^2(\Omega)} \le c_2 < \infty$ and such that,
if $x_1, x_2, \ldots$ is any sequence in $B$, then we have
$$\sum_{k=1}^n x_k \zeta_k \ \ \mbox{converges in} \ \ L^2(\Omega)
\Longrightarrow \sum_{k=1}^{\infty} \|x_k\|_B^2 < \infty.$$
\end{itemize}
\end{lemma}

\begin{lemma} \label{Quantized-Classic}
The following assertions hold:
\begin{itemize}
\item[$1.$] Let $E$ be an operator space having $\mathcal{R}$-type
$2$, then the underlying Banach space has Rademacher type $2$.
\item[$2.$] Let $E$ be an operator space having $\mathcal{R}$-cotype
$2$, then the underlying Banach space has Rademacher cotype $2$.
\end{itemize}
\end{lemma}

\emph{Proof}. Let us take a countable subset $\{\sig_k: k \ge 1\}$
of $\Sig$. Then we define the random variables $\zeta_k =
\sqrt{d_{\sig_k}} \rho_{11}^{\sig_k}$ for $k \ge 1$. The sequence
$\zeta_1, \zeta_2, \ldots$ is orthonormal in $L^2(\Omega)$ and is
made up of mean zero independent random variables. Moreover, if we
take any square-summable sequence $x_1, x_2, \ldots$ in $E$ and
$A^k \in M_{d_{\sig_k}} \otimes E$ is defined by $A_{ij}^k =
\delta_{i1} \delta_{j1} \,\ d_{\sig_k}^{-1/2} x_k$, we have
\begin{eqnarray*}
\Big\| \sum_{k=m_1+1}^{m_2} x_k \zeta_k \Big\|_{L_E^2(\Omega)} & =
& \Big\| \sum_{k=m_1+1}^{m_2} d_{\sig_k} \mbox{tr} (A^k
\rho^{\sig_k}) \Big\|_{L_E^2(\Omega)} \\ & \le &
\mathcal{K}_{12}(E, \mathcal{R}) \,\ \Big( \sum_{k=m_1+1}^{m_2}
d_{\sig_k} \|A^k\|_{S_{d_{\sig_k}}^2(E)}^2 \Big)^{1/2} \\ & = &
\mathcal{K}_{12}(E, \mathcal{R}) \,\ \Big( \sum_{k=m_1+1}^{m_2}
\|x_k\|_E^2 \Big)^{1/2}.
\end{eqnarray*}
Similarly, we get $$\Big( \sum_{k=m_1+1}^{m_2} \|x_k\|_E^2
\Big)^{1/2} \le \,\ \mathcal{K}_{22}(E, \mathcal{R}) \,\ \Big\|
\sum_{k=m_1+1}^{m_2} x_k \zeta_k \Big\|_{L_E^2(\Omega)}.$$ That
is, we have proved that $$\sum_{k=1}^{\infty} \|x_k\|_B^2 < \infty
\Longleftrightarrow \sum_{k=1}^n x_k \zeta_k \ \ \mbox{converges
in} \ \ L_E^2(\Omega).$$ But convergence in $L_E^2(\Omega)$
implies a.s. convergence for these kind of series, see Theorem
2.10 in Chapter 3 of \cite{AG}. The proof is concluded by applying
lemma \ref{Araujo-Gine}. $\blacksquare$

\begin{remark} \emph{By duality, a similar result holds for
the strong cotype.}
\end{remark}

In this section we explore Kwapie\'n theorem for the present
context. That is, completely isomorphic characterizations of
Pisier's OH Hilbertian operator spaces by means of quantized
orthonormal systems. Roughly speaking, an OH operator space is the
only possible quantization on a Hilbert space such that the
canonical identification between the resulting operator space and
its antidual is a complete isometry, see \cite{P1} for a complete
study of these spaces. In other words, the OH operator spaces are
the natural substitutes of classical Hilbert spaces in the
category of operator spaces.

\begin{theorem} \label{Kwapien1}
Let $\Phi$ be any u.b.q.o.s. associated to the parameters $(\Sig,
\mathbf{d}_{\Sig})$. Let $E$ be an operator space, then the
following are equivalent:
\begin{itemize}
\item[$1.$] $E$ is completely isomorphic to some \textnormal{OH}
Hilbertian operator space.
\item[$2.$] $E$ has $\Phi$-type and $\Phi$-cotype $2$.
\end{itemize}
\end{theorem}

\emph{Proof}. We begin by showing $(1 \Rightarrow 2)$. Let us
assume that $E$ is completely isomorphic to OH(I), for some index
set I. Then we invoke the general results stated in section
\ref{Section-Definitions} to write $\mathcal{K}_{12}(E,\Phi) \le
d_{cb}(E, \mbox{OH(I)}) \,\ \mathcal{K}_{12}(\mbox{OH(I)}, \Phi)$.
But OH(I) is completely isometric to $l^2(I)$ and it is not
difficult to check that $\mathcal{K}_{12}(l^2(I), \Phi) = 1$. This
shows that $E$ has $\Phi$-type $2$. Similar arguments are valid to
see that $E$ also has $\Phi$-cotype $2$.

Now we see $(2 \Rightarrow 1)$. Let us suppose that $E$ has
$\Phi$-type and $\Phi$-cotype $2$. By proposition \ref{Extremal}
we can replace $\Phi$ by the quantized Rademacher system
$\mathcal{R}$ of parameters $(\Sig, \mathbf{d}_{\Sig})$. We know
that $S^2(E)$ also has $\mathcal{R}$-type and $\mathcal{R}$-cotype
$2$, see again the general results of section
\ref{Section-Definitions}. Now, lemma \ref{Quantized-Classic}
gives that (the underlying Banach space of) $S^2(E)$ has
Rademacher type and cotype $2$. In particular, $S^2(E)$ is
isomorphic to some Hilbert space. By Kwapie\'n's original theorem,
this geometric condition on $S^2(E)$ is equivalent to the
existence of a constant $c$ such that $\|T \otimes
I_{S^2(E)}\|_{\mathcal{B}(l_n^2(S^2(E)), \hspace{2pt}
l_n^2(S^2(E)))} \le c \,\ \|T\|_{\mathcal{B}(l_n^2, \hspace{2pt}
l_n^2)}$ for any linear mapping $T: l_n^2 \rightarrow l_n^2$ and
any $n \ge 1$, see \cite{K1}. On the other hand, the Fubini
complete isometry $l_n^2(S^2(E)) \simeq S^2(l_n^2(E))$ given in
\cite{P2} allows us to write the last inequality as $\|T \otimes
I_E\|_{cb(l_n^2(E), \hspace{2pt} l_n^2(E))} \le c \,\
\|T\|_{\mathcal{B}(l_n^2, \hspace{2pt} l_n^2)}$. Finally, Pisier
proved that this condition is equivalent to condition $1$, see
Theorem $6.11$ of \cite{P2}. This completes the proof.
$\blacksquare$

\section{Complete quantized orthonormal systems}
\label{Section-Kwapien2}

In this section we extend the operator space version of Kwapie\'n
theorem to complete quantized orthonormal systems, uniformly
bounded or not. To that end we begin by recalling that, since
$(\Omega, \mathcal{M}, \mu)$ has no atoms, we can define a family
of dyadic sets $\mathrm{D}_j^k$ on $\Omega$, where $1 \le j \le
2^k$ and $k \ge 1$, satisfying the following conditions:
\begin{itemize}
\item $\mathrm{D}_j^k = \mathrm{D}_{2j-1}^{k+1} \cup
\mathrm{D}_{2j}^{k+1}$ for all $k \ge 1$ and all $1 \le j \le
2^k$.
\item $\Omega$ is the disjoint union of $\mathrm{D}_j^k$ for any fixed
$k \ge 1$ and all $1 \le j \le 2^k$.
\item The sets $\mathrm{D}_j^k$ are $\mu$-measurable and
$\mu(\mathrm{D}_j^k) = 2^{-k}$.
\end{itemize}
Then, if $\mbox{{\large 1}}_{\Lambda}$ stands for the
characteristic function of a measurable set $\Lambda \subset
\Omega$, we define the system $\Delta$ on $L^2(\Omega)$ by the
functions $$\delta_k = \sum_{j=1}^{2^k} (-1)^{j+1} \mbox{{\large
1}}_{\mathrm{D}_j^k}.$$

\begin{lemma} \label{Approximation}
Let $\Psi = \{\psi^{\sig}: \Omega \rightarrow M_{\Degree}\}_{\sig
\in \Sig}$ be a complete quantized orthonormal system. Let
$\{\varepsilon_n: n \ge 1\}$ be any sequence of positive numbers.
Then there exists a sequence $f_1, f_2, \ldots$ of
$\Psi$-polynomials and an increasing subsequence $k_1, k_2 \ldots$
of positive integers satisfying:
\begin{itemize}
\item[$1.$] $\Fourier_{\Psi}(f_1), \Fourier_{\Psi}(f_2), \ldots$ have
pairwise disjoint supports on $\Sig$.
\item[$2.$] $\|f_n - \delta_{k_n}\|_{L^2(\Omega)} < \varepsilon_n$.
\end{itemize}
\end{lemma}

\emph{Proof}. Let $\sig \in \Sig$ and let us fix $1 \le i,j \le
\Degree$. Then, since $\Delta$ is orthonormal in $L^2(\Omega)$,
Bessel inequality provides the following estimate
$$\sum_{k=1}^{\infty} \big| \Fourier_{\Psi}(\delta_k)_{ij}^{\sig}
\big|^2 = \sum_{k=1}^{\infty} \Big| \int_{\Omega} \delta_k
(\omega) \overline{\psi_{ji}^{\sig}(\omega)} d \mu(\omega) \Big|^2
\le \|\psi_{ji}^{\sig}\|_{L^2(\Omega)}^2 = \frac{1}{\Degree} <
\infty.$$ In particular, for all $\epsilon > 0$ and for all finite
subset $\Sig_0 \subset \Sig$ there exists a positive integer
$m(\Sig_0,\epsilon)$ such that for all $k \ge m(\Sig_0,\epsilon)$
we have $$\sum_{\sig \in \Sig_0} \Degree \sum_{i,j = 1}^{\Degree}
\big| \Fourier_{\Psi}(\delta_k)_{ij}^{\sig} \big|^2 < \epsilon.$$
On the other hand, let $\Psi_0$ be the space of
$\Psi$-polynomials. That is, $\Psi_0$ is the span of the entries
$\psi_{ij}^{\sig}$ where $1 \le i,j \le \Degree$ and $\sig \in
\Sig$. Then we construct the functions $f_1, f_2, \ldots$ as
follows:
\begin{itemize}
\item Let $f_1 \in \Psi_0$ be such that
$\|f_1 - \delta_1\|_{L^2(\Omega)} < \varepsilon_1$.
\item For $n > 1$, let $\epsilon_n = \varepsilon_n / 3$ and let
$$\Sig_n = \bigcup_{k=1}^{n-1} \,\ \mbox{supp} \big(
\Fourier_{\Psi}(f_k) \big) \subset \Sig.$$ If $k_n = m(\Sig_n,
\epsilon_n)$ we take $g_n$ to be any function in $\Psi_0$
satisfying the estimate $\|g_n - \delta_{k_n}\|_{L^2(\Omega)} <
\epsilon_n$. Then we define $$f_n = g_n - \sum_{\sig \in \Sig_n}
\Degree \mbox{tr} \big( \Fourier_{\Psi}(g_n)^{\sig} \,\
\psi^{\sig} \big).$$
\end{itemize}
The verification that the sequence $f_1, f_2, \ldots$ satisfies
the required properties is left to the reader. This completes the
proof. $\blacksquare$

\begin{theorem} \label{Kwapien2}
Let $\Psi$ be any complete quantized orthonormal system with
parameters $(\Sig, \mathbf{d}_{\Sig})$. Let $E$ be an operator
space, the following are equivalent:
\begin{itemize}
\item[$1.$] $E$ is completely isomorphic to some \textnormal{OH}
Hilbertian operator space.
\item[$2.$] $E$ has $\Psi$-type and $\Psi$-cotype $2$.
\end{itemize}
\end{theorem}

\emph{Proof}. The arguments used in theorem \ref{Kwapien1} to see
$(1 \Rightarrow 2)$ are also valid here. Let us prove that $(2
\Rightarrow 1)$. First we recall that, by lemma
\ref{Approximation}, there exists a sequence $f_1, f_2, \ldots$ of
$\Psi$-polynomials $$f_n = \sum_{\sig \in \Sig_n} \sum_{1 \le i,j
\le \Degree} \alpha_{ij}^{\sig} \psi_{ij}^{\sig}$$ where
$\alpha_{ij}^{\sig} \in \C$, $\Sig_n$ is some finite subset of
$\Sig$ and such that
\begin{itemize}
\item $\Sig_{n_1} \cap \Sig_{n_2} = \emptyset \,\ $ whenever $n_1 \neq n_2$.
\item $\|f_n - \delta_{k_n}\|_{L^2(\Omega)} < 2^{-n}$.
\end{itemize}
Now, if $E$ has $\Psi$-type and $\Psi$-cotype $2$ then the same
holds for $F = S^2(E)$. In particular, for any family $\{x_1, x_2,
\ldots x_n\}$ in $F$, we have $$\Big\| \sum_{n=1}^m x_n
\delta_{k_n} \Big\|_{L_F^2(\Omega)} \le \Big\| \sum_{n=1}^m x_n
(\delta_{k_n} - f_n) \Big\|_{L_F^2(\Omega)} + \Big\| \sum_{n=1}^m
x_n f_n \Big\|_{L_F^2(\Omega)} = \mathrm{A} + \mathrm{B}.$$ By
H\"{o}lder's inequality we get $$\mathrm{A} \le \Big( \sum_{n=1}^m
\|f_n - \delta_{k_n}\|_{L^2(\Omega)}^2 \Big)^{1/2} \Big(
\sum_{n=1}^m \|x_n\|_F^2 \Big)^{1/2} \le \frac{1}{\sqrt{3}} \Big(
\sum_{n=1}^m \|x_n\|_F^2 \Big)^{1/2}.$$ And, in order to estimate
$\mathrm{B}$, we write
\begin{eqnarray*}
\mathrm{B} & = & \Big\| \sum_{n=1}^m x_n \sum_{\sig \in \Sig_n}
\sum_{1 \le i,j \le \Degree} \alpha_{ij}^{\sig} \,\
\psi_{ij}^{\sig} \Big\|_{L_F^2(\Omega)} \\ & = & \Big\|
\sum_{n=1}^m \sum_{\sig \in \Sig_n} \Degree \,\ \mbox{tr} \big[
(\Fourier_{\Psi}(f_n)^{\sig} \otimes x_n ) \,\ \psi^{\sig} \big]
\Big\|_{L_F^2(\Omega)} \\ & \le & \mathcal{K}_{12}(E, \Psi) \,\
\Big( \sum_{n=1}^m \|x_n\|_F^2 \sum_{\sig \in \Sig_n} \Degree
\big\| \Fourier_{\Psi}(f_n)^{\sig} \big\|_{S_{\Degree}^2}^2
\Big)^{1/2}
\\ & = & \mathcal{K}_{12}(E, \Psi) \,\ \Big( \sum_{n=1}^m
\|x_n\|_F^2 \,\ \|f_n\|_{L^2(\Omega)}^2 \Big)^{1/2}
\\ & \le & 2 \,\ \mathcal{K}_{12}(E, \Psi) \,\ \Big( \sum_{n=1}^m
\|x_n\|_F^2 \Big)^{1/2}.
\end{eqnarray*}
That is, if $\Delta'$ stands for the system in $L^2(\Omega)$
defined by the functions $\delta_{k_1}, \delta_{k_2}, \ldots$,
then we have shown that $F$ has $\Delta'$-type $2$ in the sense of
\cite{GKKT}. But this is equivalent to saying that $E$ has
$\Delta'$-type $2$ in the sense of definition \ref{Riesz}. Similar
arguments are valid to see that $E$ also has $\Delta'$-cotype $2$.
Then the proof is concluded by applying theorem \ref{Kwapien1}.
$\blacksquare$

\begin{remark} \label{Failure}
\emph{The analog of Kwapie\'n's argument given in \cite{K1} for
this result does not work. Namely, if $\mathcal{R}$ denotes the
quantized Rademacher system with parameters $(\Sig,
\mathbf{d}_{\Sig})$, the idea is to use the completeness of $\Psi$
to construct a sequence $f^{\sig_1}, f^{\sig_2}, \ldots$ of
matrix-valued $\Psi$-polynomials with non-overlapping ranges of
frequencies and such that $$\int_{\Omega} \|\rho^{\sig_n} -
f^{\sig_n}\|_{S_{d_{\sig_n}}^2}^2 d \mu < \varepsilon_n \ \
\mbox{with} \ \ \varepsilon_1, \varepsilon_2, \ldots \ \
\mbox{small enough}.$$ This sequence exists and its construction
is similar to the one provided in lemma \ref{Approximation}. If
$\mathcal{R}'$ denotes the subsystem of $\mathcal{R}$ defined by
the functions $\rho^{\sig_1}, \rho^{\sig_2}, \ldots$, the next
step is to show that $\Psi$-type $2$ implies $\mathcal{R}'$-type
$2$ and the same for the cotype. Here is where the proof fails.
However, it can be checked that it works in the following cases:
\begin{itemize}
\item $\Psi$-type $2$ $\Rightarrow$ $\mathcal{R}'$-type $2$
if $\mathbf{d}_{\Sig}$ is bounded.
\item $\Psi$-cotype $2$ $\Rightarrow$ $\mathcal{R}'$-cotype $2$ if
$\Degree = 1$ for all $\sig \in \Sig$.
\end{itemize}}
\end{remark}

\section{The probabilistic approach}
\label{Section-Kwapien3}

In this section we introduce the quantization of the classical
Gauss system and analyze its important role in the operator space
version of Kwapie\'n theorem. First we outline a simple proof of
Kwapie\'n theorem for this system and then we give an alternative
proof following Kwapie\'n's approach in \cite{K1} conveniently
adapted to our setting. The reason for this approach will be clear
in corollary \ref{Kwapien-Banach}.

\begin{definition}
\emph{ Let $\{\gamma_{ij}^{\sig}: \Omega \rightarrow \mathbb{R},
\,\ 1 \le i,j \le \Degree\}_{\sig \in \Sig}$ be a family of
independent real gaussian random variables with mean zero and
variance $1$. Then the collection $\mathcal{G} = \{\gamma^{\sig}:
\Omega \rightarrow M_{\Degree}\}_{\sig \in \Sig}$, where
$\gamma^{\sig}$ stands for the random matrix $$\gamma^{\sig} =
\frac{1}{\sqrt{\Degree}} \,\ \Big(\,\ \gamma_{ij}^{\sig} \,\
\Big)$$ defines the quantized \emph{gaussian} system associated to
$(\Sig, \mathbf{d}_{\Sig})$.}
\end{definition}

\begin{remark}
\emph{Analogously, considering a priori complex gaussian random
variables, we get the quantized \emph{complex gaussian} system
associated to $(\Sig, \mathbf{d}_{\Sig})$.}
\end{remark}

This quantized system satisfies orthonormality but fail to be
uniformly bounded or complete. So all the previous results do not
seem to be valid for the quantized gaussian system. However, it is
not difficult to check that lemma \ref{Quantized-Classic} remains
valid when we replace the quantized Rademacher system
$\mathcal{R}$ by the quantized gaussian system $\mathcal{G}$. In
particular, the proof of theorem \ref{Kwapien1} also holds for
$\mathcal{G}$.

We are giving an alternative approach to this result. Let
$\widetilde{\Omega}$ be the proba\-bility space formed by the
product of infinitely many copies of $\Omega$ $$\widetilde{\Omega}
= \prod_{k=1}^{\infty} \,\ \Omega_k \qquad \mbox{and} \qquad
\widetilde{\mu} = \prod_{k=1}^{\infty} \,\ \mu_k$$ with $\Omega_k
= \Omega$ and $\mu_k = \mu$ for all $k \ge 1$. The random matrix
$\rho^{\sig,k}: \widetilde{\Omega}: \rightarrow O(d_{\sig})$ is
defined as a copy of $\rho^{\sig}$, the $\sig$-th Rademacher
function, depending only on the $k$-th coordinate. Also, for each
positive integer $m$, we define $$\rho^{\sig}(m):
\widetilde{\Omega} \longrightarrow M_{\Degree} \quad \mbox{as}
\quad \rho^{\sig}(m) = \frac{1}{\sqrt{m}} \sum_{k=1}^m
\rho^{\sig,k}.$$ Finally, we construct a quantized gaussian system
$\{\widetilde{\gamma}^{\sig}: \widetilde{\Omega} \rightarrow
M_{\Degree}\}_{\sig \in \Sig}$ on $\widetilde{\Omega}$ associated
to the parameters $(\Sig, \mathbf{d}_{\Sig})$. We state a slight
modification of the central limit theorem in type $2$ spaces, see
\cite{AG} for the classical statement of that result. It is
nothing but an analog, for Banach-valued random variables, of
Lemma $2.1$ in \cite{K1}. Let us fix a finite subset $\Sig_0 =
\{\sig_1, \sig_2, \ldots \sig_n\}$ of $\Sig$.

\begin{proposition} \label{Central-Limit}
Let $\displaystyle h: S_{d_{\sig_1}}^2 \times \cdots \times
S_{d_{\sig_n}}^2 \rightarrow \mathbb{R}$ be a continuous function
such that
\begin{equation} \label{Control}
h(D^{\sig_1}, \ldots, D^{\sig_n}) \,\ e^{\mbox{{\scriptsize
$\displaystyle - \sum_{j=1}^n \|D^{\sig_j}\|$}}} \longrightarrow 0
\qquad \mbox{as} \quad \sum_{j=1}^n
\|D^{\sig_j}\|_{S_{d_{\sig_j}}^2} \rightarrow \infty.
\end{equation}
Then we have $\displaystyle \ \ \lim_{m \rightarrow \infty}
\int_{\widetilde{\Omega}} h(\rho^{\sig_1}(m), \ldots,
\rho^{\sig_n}(m))) d \widetilde{\mu} = \int_{\widetilde{\Omega}}
h(\widetilde{\gamma}^{\sig_1}, \ldots,
\widetilde{\gamma}^{\sig_n}) d \widetilde{\mu}$.
\end{proposition}

\emph{Sketch of the proof}. By using the orthonormality relations
of quantized Rademacher and gaussian systems, one easily gets that
the distribution of $\widetilde{\gamma}^{\sig}$ is a centered
cylindrical gaussian measure with the same covariance as that of
$\rho^{\sig,k}$ for all $k \ge 1$. Hence, by the central limit
theorem in type $2$ spaces, the joint distribution of
$(\rho^{\sig_1}(m), \ldots, \rho^{\sig_n}(m))$ converges weakly to
the joint distribution of $(\widetilde{\gamma}^{\sig_1}, \ldots,
\widetilde{\gamma}^{\sig_n})$. Now, if we write $S_{\Sig_0}^2 =
S_{d_{\sig_1}}^2 \times \cdots \times S_{d_{\sig_n}}^2$, we define
the Banach space $B$ of all continuous functions $h: S_{\Sig_0}^2
\rightarrow \mathbb{R}$ satisfying (\ref{Control}) and with the
norm given by $$\|h\|_{B} = \sup \Big\{ |h(D^{\sig_1}, \ldots,
D^{\sig_n})| \,\ e^{\mbox{{\scriptsize $\displaystyle -
\sum_{j=1}^n \|D^{\sig_j}\|$}}}: \,\ (D^{\sig_1}, \ldots,
D^{\sig_n}) \in S_{\Sig_0}^2 \Big\}.$$ We also define the
following functionals on $B$
\begin{eqnarray*}
T(h) & = & \int_{\widetilde{\Omega}}
h(\widetilde{\gamma}^{\sig_1}, \ldots,
\widetilde{\gamma}^{\sig_n}) d \widetilde{\mu} \ \ \mbox{and} \ \
T_m(h) = \int_{\widetilde{\Omega}} h(\rho^{\sig_1}(m), \ldots,
\rho^{\sig_n}(m))) d \widetilde{\mu}.
\end{eqnarray*}
Following the arguments given in Lemma $2.1$ of \cite{K1}, it
suffices to check that $T$ and $T_m$ are well-defined and that
$\sup \|T_m\| < \infty$. $T_m$ is well-defined since
$h(\rho^{\sig_1}(m), \ldots, \rho^{\sig_n}(m))$ is a bounded
function. On the other hand,
\begin{eqnarray*}
|T(h)| & \le & \|h\|_B \,\ \prod_{j=1}^n \,\
\int_{S_{d_{\sig_j}}^2} \exp \|D^{\sig_j}\|_{S_{d_{\sig_j}}^2} \,\
d \mu_{\widetilde{\gamma}^{\sig_j}} (D^{\sig_j}) \\ & \le &
\|h\|_B \,\ \prod_{j=1}^n \,\ \int_{\widetilde{\Omega}} \prod_{1
\le i_1, i_2 \le d_{\sig_j}} \exp \Big|
\frac{\widetilde{\gamma}_{i_1i_2}^{\sig_j}(\widetilde{\omega})}{\sqrt{d_{\sig_j}}}
\Big| \,\ d \widetilde{\mu}(\widetilde{\omega}) \\ & \le & \|h\|_B
\,\ \prod_{j=1}^n \prod_{1 \le i_1, i_2 \le d_{\sig_j}} \Big(
\int_{\widetilde{\Omega}} \exp d_{\sig_j}^2 \Big|
\frac{\widetilde{\gamma}_{i_1i_2}^{\sig_j}(\widetilde{\omega})}{\sqrt{d_{\sig_j}}}
\Big| \,\ d \widetilde{\mu}(\widetilde{\omega})
\Big)^{1/d_{\sig_j}^2} \\ & = & \|h\|_B \,\ \prod_{j=1}^n \prod_{1
\le i_1, i_2 \le d_{\sig_j}} \Big( \frac{1}{\sqrt{2 \pi}}
\int_{\mathbb{R}} \exp |d_{\sig_j}^{3/2} s| \exp(-s^2/2) ds
\Big)^{1/d_{\sig_j}^2}
\end{eqnarray*}
where we have applied the obvious inequality $\|D\|_{S_n^2} \le
\sum_{ij} |D_{ij}|$ and the generalized H\"{o}lder inequality.
Therefore $T$ is well-defined. Similar arguments give the uniform
boundedness of $\|T_m\|$. $\blacksquare$

\

Before the proof of Kwapie\'n theorem for the quantized gaussian
system we need to state a couple of lemmas. Let $D_1$, $D_2$ be
orthogonal $\Degree \times \Degree$ matrices, then $D_1
\gamma^{\sig} D_2$ and $\gamma^{\sig}$ have the same distribution.
The next result can be found in \cite{MP}, it follows from this
\lq sign invariance\rq ${}$ and the contraction principle stated
above.

\begin{lemma} \label{Marcus-Pisier}
Let $B$ be a Banach space. There exists a positive constant $c$,
such that for any finite set $\Sig_0$ of $\Sig$, we have
$$\int_{\Omega} \Big\| \sum_{\sig \in \Sig_0} \Degree
\textnormal{tr} (A^{\sig} \rho^{\sig}) \Big\|_B^2 d \mu \le c
\int_{\Omega} \Big\| \sum_{\sig \in \Sig_0} \Degree
\textnormal{tr} (A^{\sig} \gamma^{\sig}) \Big\|_B^2 d \mu.$$
\end{lemma}

The following result is a completely isomorphic characterization
of OH operator spaces given by Pisier in \cite{P2}. It can be
regarded as the version for operator spaces of a previous
isomorphic characterization of Hilbert spaces given by Kwapie\'n,
see $(\textrm{iv})$ of Proposition $3.1$ in \cite{K1}.

\begin{lemma} \label{Pisier}
Let $E$ be an operator space. Then $E$ is completely isomorphic to
some \textnormal{OH} Hilbertian operator space if and only if
there exists a positive constant $c$ such that for any $n \ge 1$
and any linear mapping $T: S_n^2 \rightarrow S_n^2$, we have $$\|T
\otimes I_E\|_{\mathcal{B}(S_n^2(E), S_n^2(E))} \le c \,\
\|T\|_{\mathcal{B}(S_n^2, S_n^2)}.$$
\end{lemma}

In the next result we assume that the gaussian system we work with
takes values in arbitrary large matrices. We need to require that
in view of the proof we are giving. Although, as we have seen,
this requirement is not necessary, it will become very natural in
corollary \ref{Kwapien-Banach}.

\begin{theorem} \label{Kwapien3}
Let $\mathcal{G}$ be the gaussian system with
parameters $(\Sig, \mathbf{d}_{\Sig})$. Let us assume that
$\mathbf{d}_{\Sig}$ is unbounded, then the following are
equivalent:
\begin{itemize}
\item[$1.$] $E$ is completely isomorphic to some \textnormal{OH}
Hilbertian operator space.
\item[$2.$] $E$ has $\mathcal{G}$-type and $\mathcal{G}$-cotype $2$.
\end{itemize}
\end{theorem}

\emph{Proof}. Let us prove that ($1 \Rightarrow 2$). Let us assume
that $E$ is completely isomorphic to some OH(I). If $\mathcal{R}$
denotes the quantized Rademacher system with parameters $(\Sig,
\mathbf{d}_{\Sig})$, then we know by theorem \ref{Kwapien1} that
$E$ has $\mathcal{R}$-type and $\mathcal{R}$-cotype $2$. But then
lemma \ref{Marcus-Pisier} gives that $E$ has $\mathcal{G}$-cotype
$2$. Let us see that $E$ also has $\mathcal{G}$-type $2$. Here we
recall that any Banach space $B$ with Rademacher type $2$
satisfies the inequality $$\int_{\Omega} \Big\| \sum_{k=1}^n
\phi_k \Big\|_B^2 d \mu \le c \sum_{k=1}^n \int_{\Omega}
\|\phi_k\|_B^2 d \mu$$ for some universal constant $c$ and any
family $\phi_1, \phi_2, \ldots \phi_n$ of mean zero independent
$B$-valued random variables in $L^2(\Omega)$. In particular, since
(by lemma \ref{Quantized-Classic}) the underlying Banach space of
$S_n^2(E)$ has Rademacher type $2$ for any $n \ge 1$, we have
\begin{eqnarray*}
\Big\| \sum_{\sig \in \Sig_0} \Degree \mbox{tr} (A^{\sig}
\rho^{\sig}(m)) \Big\|_{S_n^2(L_E^2(\widetilde{\Omega}))} & = &
\Big\| \sum_{k=1}^m \sum_{\sig \in \Sig_0} \Degree \mbox{tr} \big[
\frac{A^{\sig} \rho^{\sig,k}}{\sqrt{m}} \big]
\Big\|_{L_{S_n^2(E)}^2(\widetilde{\Omega})} \\ & \le & c \,\ \Big(
\sum_{k=1}^m \Big\| \sum_{\sig \in \Sig_0} \Degree \mbox{tr} \big[
\frac{A^{\sig} \rho^{\sig,k}}{\sqrt{m}} \big]
\Big\|^2_{L_{S_n^2(E)}^2(\Omega_k)} \Big)^{1/2} \\ & \le & c \,\
\mathcal{K}_{12}(E, \mathcal{R}) \,\ \Big( \sum_{\sig \in \Sig_0}
\Degree \|A^{\sig}\|_{S_{\Degree n}^2}^2 \Big)^{1/2}
\end{eqnarray*}
On the other hand, let $h(D^{\sig_1}, \ldots, D^{\sig_n}) =
\|\sum_{\sig \in \Sig_0} \Degree \mbox{tr}(A^{\sig}
D^{\sig})\|_{S_n^2(E)}^2$. Let us see that $h$ satisfies
hypothesis (\ref{Control}) of proposition \ref{Central-Limit}.
First we recall that
\begin{eqnarray*}
\|\mbox{tr}(A^{\sig}D^{\sig})\|_{S_n^2(E)} & = &
\sup_{\|T\|_{S_n^2(E^{\star})} \le 1} \,\ \mbox{tr} \big( A^{\sig}
(D^{\sig} \otimes T) \big) \\ & \le &
\sup_{\|T\|_{S_n^2(E^{\star})} \le 1} \,\ \|A^{\sig}\|_{S_{\Degree
n}^2(E)} \|D^{\sig} \otimes T\|_{S_{\Degree n}^2(E^{\star})} \\ &
= & \|A^{\sig}\|_{S_{\Degree n}^2(E)}
\|D^{\sig}\|_{S_{\Degree}^2}.
\end{eqnarray*}
Hence we get
\begin{eqnarray*} h(D^{\sig_1}, \ldots , D^{\sig_n}) & \le
& \Big( \sum_{\sig \in \Sig_0} \Degree \|A^{\sig}\|_{S_{\Degree
n}^2(E)} \|D^{\sig}\|_{S_{\Degree}^2} \Big)^2 \\ & \le &
\max_{\sig \in \Sig_0} \Degree^2 \|A^{\sig}\|_{S_{\Degree
n}^2(E)}^2 \,\ \Big( \sum_{\sig \in \Sig_0}
\|D^{\sig}\|_{S_{\Degree}^2} \Big)^2
\end{eqnarray*}
and so $h$ satisfies (\ref{Control}). In particular, we apply
proposition \ref{Central-Limit} to obtain
\begin{eqnarray*}
\Big\| \sum_{\sig \in \Sig_0} \Degree \mbox{tr}(A^{\sig}
\gamma^{\sig}) \Big\|_{S_n^2(L_{E}^2(\Omega))} & = & \Big\|
\sum_{\sig \in \Sig_0} \Degree \mbox{tr}(A^{\sig}
\widetilde{\gamma}^{\sig})
\Big\|_{S_n^2(L_E^2(\widetilde{\Omega}))} \\ & = & \lim_{m
\rightarrow \infty} \Big\| \sum_{\sig \in \Sig_0} \Degree
\mbox{tr}(A^{\sig} \rho^{\sig}(m))
\Big\|_{S_n^2(L_E^2(\widetilde{\Omega}))} \\ & \le & c \,\
\mathcal{K}_{12}(E, \mathcal{R}) \,\ \Big( \sum_{\sig \in \Sig_0}
\Degree \|A^{\sig}\|_{S_{\Degree n}^2}^2 \Big)^{1/2}
\end{eqnarray*}
Therefore, by Lemma $1.7$ of \cite{P2}, we obtain that $E$ has
$\mathcal{G}$-type $2$ and the proof of $(1 \Rightarrow 2)$ in
concluded.

Now we see ($2 \Rightarrow 1$). By the unboundedness of
$\mathbf{d}_{\Sig}$ and lemma \ref{Pisier} it suffices to see that
there exists a positive constant $c$ such that, for any $\sig \in
\Sig$ and any linear mapping $T: S_{\Degree}^2 \rightarrow
S_{\Degree}^2$, we have
\begin{equation} \label{Condition}
\|T \otimes I_{E}\|_{\mathcal{B}(S_{\Degree}^2(E),
S_{\Degree}^2(E))} \le c \,\ \|T\|_{\mathcal{B}(S_{\Degree}^2,
S_{\Degree}^2)}.
\end{equation}
By homogeneity it is enough to see (\ref{Condition}) for $T$ in
the unit ball $\mathrm{B}_{\sig}$ of $\mathcal{B}(S_{\Degree}^2,
S_{\Degree}^2)$. But $\mathrm{B}_{\sig}$ is a compact, convex set
and then every element of $\mathrm{B}_{\sig}$ is a convex linear
combination of unitary operators, the extreme points of
$\mathrm{B}_{\sig}$. Therefore, it suffices to check
(\ref{Condition}) for $T$ unitary. Let $A \in S_{\Degree}^2(E)$
and $T: S_{\Degree}^2 \rightarrow S_{\Degree}^2$ unitary, then we
have
\begin{eqnarray*}
\|T \otimes I_E(A)\|_{S_{\Degree}^2(E)} & \le & \Degree^{-1/2} \,\
\mathcal{K}_{22}(E, \mathcal{G}) \,\ \|\Degree \hspace{2pt}
\mbox{tr}(\gamma^{\sig} \hspace{2pt} [T \otimes I_E]
(A))\|_{L_E^2(\Omega)}
\\ & = & \Degree^{-1/2} \,\ \mathcal{K}_{22}(E, \mathcal{G}) \,\
\|\Degree \hspace{2pt} \mbox{tr}(T^{\star}(\gamma^{\sig})
\hspace{2pt} A)\|_{L_E^2(\Omega)}
\\ & = & \Degree^{-1/2} \,\ \mathcal{K}_{22}(E, \mathcal{G}) \,\
\|\Degree \hspace{2pt} \mbox{tr}(\gamma^{\sig} \hspace{2pt}
A)\|_{L_E^2(\Omega)}
\\& \le & \mathcal{K}_{22}(E, \mathcal{G}) \,\ \mathcal{K}_{12}(E,
\mathcal{G}) \,\ \|A\|_{S_{\Degree}^2(E)}
\end{eqnarray*}
since, by the unitarity of $T$, the distribution of
$T(\gamma^{\sig})$ is the same as that of $\gamma^{\sig}$ (see
Theorem $6.8$ in Chapter $3$ of \cite{AG}). Therefore $E$
satisfies condition (\ref{Condition}). This completes the proof.
$\blacksquare$

\

Let $\Phi$ be a quantized orthonormal system and let $E$ be an
operator space. Let $1 \le p \le 2$, we shall say that $E$ has
\emph{Banach $\Phi$-type} $p$ if
$$\widetilde{\mathcal{K}}_{1p}(E,\Phi) = \sup
\|\Fourier_{\Phi_0}^{-1} \otimes
I_E\|_{\mathcal{B}(\Lebesgue_E^p(\Sig_0), \Phi_E^{p'}(\Sig_0))} <
\infty$$ where the supremum runs over the family of finite subsets
$\Sig_0$ of $\Sig$. That is, we do not require the complete
boundedness of $\Fourier_{\Phi_0}^{-1} \otimes I_E$ as we did in
definition \ref{Riesz}, we just require the boundednes of it. In
the same fashion can be defined the \emph{Banach $\Phi$-cotype}
$p'$ of an operator space and the subsequent constant
$\widetilde{\mathcal{K}}_{2p'}(E,\Phi)$. The following result,
which is a consequence of the probabilistic argument employed in
the proof of theorem \ref{Kwapien3}, shows that the notions of
Banach $\Phi$-type and Banach $\Phi$-cotype $2$ are the right ones
in the operator space version of Kwapie\'n's theorem whenever the
quantized system $\Phi$ takes values in arbitrary large matrices.

\begin{corollary} \label{Kwapien-Banach}
Let $\mathbf{d}_{\Sig}$ be an unbounded family of positive
integers indexed by $\Sig$. Let $\Phi$ be any u.b.q.o.s. with
parameters $(\Sig, \mathbf{d}_{\Sig})$. Let $E$ be an operator
space, then the following are equivalent:
\begin{itemize}
\item[$1.$] $E$ is completely isomorphic to some \textnormal{OH}
Hilbertian operator space.
\item[$2.$] $E$ has Banach $\Phi$-type and Banach $\Phi$-cotype
$2$.
\item[$3.$] $E$ has Banach $\mathcal{G}$-type and Banach
$\mathcal{G}$-cotype $2$.
\end{itemize}
\end{corollary}

\emph{Proof}. The implication $(1 \Rightarrow 2)$ is obvious. Now,
$(2 \Rightarrow 3)$ follows from proposition \ref{Extremal} and
the probabilistic proof of theorem \ref{Kwapien3}. Recall that the
proofs given for both results are still valid when complete
boundedness is replaced by boundedness. Finally, $(3 \Rightarrow
1)$ since the proof of theorem \ref{Kwapien3} only uses that $E$
has Banach $\mathcal{G}$-type and Banach $\mathcal{G}$-cotype $2$.
$\blacksquare$

\begin{remark}
\emph{Obviously this result fails for $\mathbf{d}_{\Sig}$ bounded.
For instance, take $\Phi$ to be the classical Rademacher system on
$L^2[0,1]$ or the dual group of the torus $\mathbb{T}$. In these
cases we go back to the classical Kwapie\'n's characterization
theorem of Hilbert spaces.}
\end{remark}

We now extend corollary \ref{Kwapien-Banach} to the case of
complete quantized orthonormal systems with $\mathbf{d}_{\Sig}$
unbounded. The proof of this result was kindly communicated to us
by Gilles Pisier.

\begin{theorem} \label{Complete}
Let $\mathbf{d}_{\Sig}$ be an unbounded family of positive
integers indexed by $\Sig$. Let $\Psi$ be any complete quantized
orthonormal system with parameters $(\Sig, \mathbf{d}_{\Sig})$.
Let $E$ be an operator space, then the following are equivalent:
\begin{itemize}
\item[$1.$] $E$ is completely isomorphic to some \textnormal{OH}
Hilbertian operator space.
\item[$2.$] $E$ has Banach $\Psi$-type and Banach $\Psi$-cotype
$2$.
\end{itemize}
\end{theorem}

\emph{Proof}. The implication $(1 \Rightarrow 2)$ is again
obvious. To see that $(2 \Rightarrow 1)$, we begin by recalling
that, if $E$ has Banach $\Psi$-type and Banach $\Psi$-cotype $2$,
then (the underlying Banach space of) $E$ is isomorphic to a
Hilbert space. That is, the proof of theorem \ref{Kwapien2} can
easily be adapted to this setting. Moreover, by another well-known
characterization of Kwapie\'n given in \cite{K2}, we know that
there exists a positive constant $c$ such that, for any linear
mapping $L: L^2(\Omega) \rightarrow L^2(\Omega)$, we have $\|L
\otimes I_E\|_{\mathcal{B}(L_E^2(\Omega), L_E^2(\Omega))} \le c
\,\ \|L\|_{\mathcal{B}(L^2(\Omega), L^2(\Omega))}$. In particular,
if $\Lambda^2$ is any closed subspace of $L^2(\Omega)$ and
$\Lambda^2(E) = \Lambda^2 \otimes E$, we get
\begin{equation} \label{Lambda}
\|L \otimes I_E\|_{\mathcal{B}(\Lambda^2(E), \Lambda^2(E))} \le c
\,\ \|L\|_{\mathcal{B}(\Lambda^2, \Lambda^2)}
\end{equation}
for any linear mapping $L: \Lambda^2 \rightarrow \Lambda^2$. Now,
for any $\sig \in \Sig$, we consider the space $\Lambda_{\sig}^2 =
\mbox{span} \{\psi_{ij}^{\sig}: \,\ 1 \le i,j \le \Degree\}$
regarded as a subspace of $L^2(\Omega)$ and the space
$\Lambda_{\sig}^2(E) = \Lambda_{\sig}^2 \otimes E$. We also need
to consider the linear isomorphism $$\begin{array}{crcl}
T_2(\sig): & S_{\Degree}^2 & \longrightarrow & \Lambda_{\sig}^2 \\
& A & \longmapsto & \Degree \mbox{tr}(A \psi^{\sig}).
\end{array}$$ The following estimates are clear
\begin{eqnarray*}
\|T_2(\sig) \otimes I_E\|_{\mathcal{B}(S_{\Degree}^2(E),
\Lambda_{\sig}^2(E))} & \le & \Degree^{1/2} \,\
\widetilde{K}_{12}(E, \Psi) \\ \|T_2(\sig)^{-1} \otimes
I_E\|_{\mathcal{B}(\Lambda_{\sig}^2(E), S_{\Degree}^2(E))} & \le &
\Degree^{-1/2} \,\ \widetilde{K}_{22}(E, \Psi).
\end{eqnarray*}
Finally, if we consider a linear mapping $T: S_{\Degree}^2
\rightarrow S_{\Degree}^2$, then we have that $T = T_2(\sig)^{-1}
\circ L_2(\sig) \circ T_2(\sig)$ where $L_2(\sig) = T_2(\sig)
\circ T \circ T_2(\sig)^{-1}$ satisfies inequality (\ref{Lambda}).
Therefore
\begin{eqnarray*}
\|T \otimes I_E\|_{\mathcal{B}(S_{\Degree}^2(E),
S_{\Degree}^2(E))} & \le & \|T_2(\sig)^{-1} \otimes I_E\| \,\
\|L_2(\sig) \otimes I_E\| \,\ \|T_2(\sig) \otimes I_E\| \\ & \le &
c \,\ \widetilde{K}_{12}(E, \Psi) \,\ \widetilde{K}_{22}(E, \Psi)
\,\ \|L_2(\sig)\|_{\mathcal{B}(\Lambda_{\sig}^2,
\Lambda_{\sig}^2)} \\ & \le & c \,\ \widetilde{K}_{12}(E, \Psi)^2
\,\ \widetilde{K}_{22}(E, \Psi)^2 \,\
\|T\|_{\mathcal{B}(S_{\Degree}^2, S_{\Degree}^2)}
\end{eqnarray*}
But then we are satisfying the hypothesis of lemma \ref{Pisier}
since $\mathbf{d}_{\Sig}$ is unbounded. This completes the proof.
$\blacksquare$

\begin{remark}
\emph{In fact, as it can be checked, the ideas behind the proof of
theorem \ref{Complete} are also valid to prove corollary
\ref{Kwapien-Banach}. In particular, the probabilistic approach
given at the beginning of this section becomes unnecessary in
order to get corollary \ref{Kwapien-Banach}. However, we have
included it since we find it as the natural source of ideas for
these results.}
\end{remark}

Let $R$ and $C$ denote the row and column operator spaces
respectively. In \cite{P1} Pisier defined natural operator space
structures on $R \cap C$ and $R + C$ in such a way that the pair
$(R \cap C, R + C)$ becomes compatible for complex interpolation.
Moreover, Pisier proved in \cite{P2} the following surprising
complete isomorphism $$(R \cap C, R+ C)_{\theta} \simeq
\mathrm{R}_p$$ where $\theta = 1 / p$ and $\mathrm{R}_p$ is, as we
defined in remark \ref{Khintchine-Kahane}, the closure in
$L^p[0,1]$ of the subspace spanned by the classical Rademacher
functions $r_1, r_2, \ldots$ endowed with its natural operator
space structure. Pisier analyzed in \cite{P2} this operator space
structure by means of the non-commutative Khintchine inequalities
previously developed by him and Lust-Piquard, see \cite{L} and
\cite{LP}. Now we use the family of operator spaces
$\{\mathrm{R}_p: 1 \le p \le \infty\}$ to illustrate some
situations:
\begin{itemize}
\item[$(\textrm{a}_1)$] Let $\Phi$ be any u.b.q.o.s. associated to
the parameters $(\Sig, \mathbf{d}_{\Sig})$ with
$\mathbf{d}_{\Sig}$ unbounded. Then $\mathrm{R}_p$ has Banach
$\Phi$-type $2$ for any $ 2 \le p < \infty$. Namely, by the
classical Khintchine inequalities the underlying Banach space of
$\mathrm{R}_p$ is isomorphic to that of $\mathrm{R}_2$ for $1 \le
p < \infty$. Moreover, the identity mapping $I: \mathrm{R}_p
\rightarrow \mathrm{R}_2$ is a complete contraction whenever $p
\ge 2$. Therefore, there exists some constant $c$ such that
\begin{eqnarray*}
\Big\| \sum_{\sig \in \Sig_0} \Degree \mbox{tr}(A^{\sig}
\varphi^{\sig}) \Big\|_{L_{\mathrm{R}_p}^2(\Omega)} & \le & c \,\
\Big\| \sum_{\sig \in \Sig_0} \Degree \mbox{tr}(A^{\sig}
\varphi^{\sig}) \Big\|_{L_{\mathrm{R}_2}^2(\Omega)} \\ & \le & c
\,\ \widetilde{\mathcal{K}}_{12}(\mathrm{R}_2, \Phi) \,\ \Big(
\sum_{\sig \in \Sig_0} \Degree
\|A^{\sig}\|_{S_{\Degree}^2(\mathrm{R}_2)}^2 \Big)^{1/2}
\\ & \le & c \,\ \widetilde{\mathcal{K}}_{12}(\mathrm{R}_2, \Phi)
\,\ \Big( \sum_{\sig \in \Sig_0} \Degree
\|A^{\sig}\|_{S_{\Degree}^2(\mathrm{R}_p)}^2 \Big)^{1/2}.
\end{eqnarray*}
Now corollary \ref{Kwapien-Banach} gives that $\mathrm{R}_p$,
although being isomorphic to a Hilbert space, can not have Banach
$\Phi$-cotype $2$ for $2 < p < \infty$ since in that cases
$\mathrm{R}_p$ is not completely isomorphic to any OH operator
space. By theorem \ref{Complete}, the same holds when we work with
any complete quantized orthonormal system $\Psi$ with
$\mathbf{d}_{\Sig}$ unbounded.
\item[$(\textrm{a}_2)$] Similarly $\mathrm{R}_p$ has Banach
$\Phi$-cotype $2$ for any $1 \le p \le 2$ but it has not Banach
$\Phi$-type $2$ unless $p=2$. By theorem \ref{Complete}, the same
holds for any complete quantized orthonormal system $\Psi$ with
$\mathbf{d}_{\Sig}$ unbounded.
\item[$(\textrm{b})$] In the commutative theory there exist some
systems for which Kwapie\'n theorem holds requiring only one of
the type $2$ or the cotype $2$ conditions. Kwapie\'n showed in
\cite{K1} that the system of characters of the torus $\mathbb{T}$
presents this kind of autoduality. Another example is given by the
system of characters of the Cantor group $\mathbb{D}$, see
\cite{GKKT} or \cite{PW} for a proof of this fact. It is easy to
see that this autoduality remains valid in our setting. For
instance, if $E$ has $\mathbb{Z}$-type $2$, then $S^2(E)$ also
does and hence it is isomorphic to some Hilbert space $H$. But
this gives that $E$ is completely isomorphic to some OH, see the
proof of theorem \ref{Kwapien1}. In particular $\mathrm{R}_p$ can
not have Fourier type $2$ or Fourier cotype $2$ with respect to
$\mathbb{T}$ or $\mathbb{D}$ unless $p=2$. On the other hand we
know that $\mathrm{R}_p$ has Banach Fourier type and Banach
Fourier cotype $2$ with respect to $\mathbb{T}$ and $\mathbb{D}$
for any $1 \le p < \infty$.

Now, it is natural to ask if there exists a non-commutative
compact group $G$ with dual object $\Gamma$ satisfying this
autoduality. That is, such that any operator space $E$ having
$\Gamma$-type $2$ or $\Gamma$-cotype $2$ is completely isomorphic
to some OH operator space. At least we know that when
$\mathbf{d}_{\Gamma}$ is unbounded, by points $(\textrm{a}_1)$ and
$(\textrm{a}_2)$, an operator space having Banach $\Gamma$-type
$2$ or Banach $\Gamma$-cotype $2$ does not have to be completely
isomorphic to any OH operator space.

At this point it also becomes natural to ask if Banach
$\Gamma$-type $2$ and $\Gamma$-type $2$ (resp. Banach
$\Gamma$-cotype $2$ and $\Gamma$-cotype $2$) are equivalent
notions as a consequence of the unboundedness of
$\mathbf{d}_{\Gamma}$. At the time of this writing, we can not
answer these questions.
\end{itemize}

\bibliographystyle{amsplain}

\begin{flushleft}
{Departamento de Matem\'{a}ticas, C-XV \\ Universidad Aut\'{o}noma de
Madrid \\ 28049 Madrid, Spain \\ E-mail: jose.garcia-cuerva@uam.es
\\ E-mail: javier.parcet@uam.es }
\end{flushleft}

\end{document}